\documentclass{article}

\usepackage{amssymb,amscd,amsmath}
\usepackage[matrix, arrow]{xy}

\newtheorem{Thm}{Theorem}
\newtheorem{Lem}[Thm]{Lemma}
\newtheorem{Prop}[Thm]{Proposition}
\newtheorem{Cor}[Thm]{Corollary}

\newenvironment{Pf}{\textbf{Proof.}\\}{\\\noindent$\blacksquare$}

\DeclareMathOperator{\GL}{GL}
\DeclareMathOperator{\Rep}{Rep}
\DeclareMathOperator{\dimv}{\underline\dim}

\DeclareMathOperator{\Ext}{Ext}
\DeclareMathOperator{\End}{End}

\DeclareMathOperator{\HA}{\mathcal H}
\DeclareMathOperator{\CA}{\mathcal C}
\DeclareMathOperator{\CM}{\mathcal{CM}}

\newcommand{\cP}{\mathcal P}
\newcommand{\cI}{\mathcal I}
\newcommand{\cR}{\mathcal R}

\newcommand{\bQ}{\mathbb Q}

\title{The Composition Monoid and the Composition Algebra at $q=0$ of the Kronecker Quiver}
\author{Andrew Hubery\\Universit\"at Bielefeld, Germany}

\begin{document}

\maketitle

\begin{abstract}
We consider the Kronecker quiver $K$ and determine the relations for the specialisation to $q=0$ of the generic composition algebra as well as those for Reineke's composition monoid. As a corollary, we deduce that the composition monoid is a proper factor of the specialisation of the composition algebra.

We also obtain a normal form for the varieties occurring in the composition monoid in terms of Schur roots.
\end{abstract}


\section{Introduction}

It has long been known that the representation theory of quivers has many close connections with the structure theory of Kac-Moody Lie algebras. The first results in this direction appeared in Kac's article \cite{Kac1}, with greater insight being afforded by the work of Ringel on the Ringel-Hall algebra of a quiver \cite{Ringel}. Green then showed that a certain natural subalgebra of the Ringel-Hall algebra, the composition algebra, is in fact isomorphic to a twisted version of the corresponding quantum group \cite{Green}.

In his recent article \cite{Rein2}, Reineke introduces the concept of the extension monoid of a quiver, whose elements are the closed, irreducible $\GL(\alpha)$-stable subsets of the representation varieties $\Rep(\alpha)$ for all $\alpha$, and whose multiplication is given by taking a `variety of extensions'. Analogously with the Ringel-Hall algebra approach, he then considers the composition monoid, which is the submonoid generated by the simple representations.

This monoid is interesting for several reasons. Its elements can be described as those representations which have a composition series with prescribed factors in a presecribed order; there are close connections with Schur roots and the notion of canonical decomposition of a dimension vector, introduced by Kac \cite{Kac2} and later studied by Schofield \cite{Scho}; finally, it is shown in \cite{Rein2} that the quantum Serre relations, specialised to $q=0$, hold in the composition monoid.

Furthermore, at the end of the his article \cite{Rein2}, Reineke poses four questions concerning the structure of the composition monoid. These relate to defining relations for the monoid and the existence of a partial normal form.

In this article we study both the composition algebra and the composition monoid in the case of the Kronecker quiver.

Using some recent results on the structure of the composition algebra for the Kronecker quiver \cite{Szanto} (see also \cite{Zhang}) we determine defining relations for the specialisation to $q=0$ of the composition algebra. In particular, we note that infinitely many such relations are needed. We also determine an alternative PBW-basis for the composition algebra which behaves particularly nicely with respect to the specialisation to $q=0$.

In the second section we concentrate on the composition monoid itself and determine its defining relations, of which there are again infinitely many. We also construct a normal form for the elements in the composition monoid, which is a special case of the partial normal form obtained in \cite{Rein2}. This can be expressed purely in terms of Schur roots and includes as a special case the canonical decomposition of a dimension vector.

In the final section, we bring together the previous results to conclude that, at least for the Kronecker quiver, the composition monoid is in fact a proper factor of the composition algebra at $q=0$. This extends the analogous result for Dynkin quivers, claimed in \cite{Rein1} and later proved in \cite{HR}, that the composition monoid is isomorphic to the specialisation of the composition algebra at $q=0$. We are also able to answer all the questions posed by Reineke.

\textit{Acknowledgements.} The author would like to thank M. Reineke and C.M. Ringel for interesting discussions concerning this article.

\section{The Composition Algebra}

Let $K$ be the Kronecker quiver $i\rightrightarrows j$ and let $k$ be an arbitrary field. For general notions concerning the representation theory of quivers, we refer the reader to the book \cite{ARS}.

We recall the following list of representatives for the isomorphism classes of indecomposable representations of $K$ over $k$:
\begin{itemize}
\item The preprojective modules $P(m)$ of dimension vector $(m,m+1)$ for $m\geq 0$;
\item The preinjective modules $I(n)$ of dimension vector $(n+1,n)$ for $n\geq 0$;
\item The regular modules $R_{mx}$ of dimension vector $(md_x,md_x)$ for $m\geq 1$ and $x\in\mathbb P^1_k$. The module $R_{mx}$ has regular length $m$ and socle the regular simple $R_x$, and $d_x$ denotes the degree of $x$, or equivalently the degree of the associated homogeneous prime ideal in $k[X,Y]$.
\end{itemize}

For $k$ a finite field with $|k|=q$, we can define the Hall algebra $\HA(kK)$ to be the $\bQ$-vector space with basis elements the isomorphism classes of representations, equipped with the product $[M][N]:=\sum_{[X]}F_{MN}^X[X]$, where
$$F_{MN}^X:=|\{Y\leq X: X/Y\cong M, Y\cong N\}|.$$

This is generally considered `too big' to study, and so we consider instead the composition algebra $\CA(kK)$. This is the subalgebra generated by the simples --- i.e. for the Kronecker quiver, the elements $[P(0)]$ and $[I(0)]$.

N.B. The algebras $\HA(kK)$ and $\CA(kK)$ are graded by dimension vector.

\begin{Thm}[Ringel, Green, Lusztig]
Define $\CA'_q(K)$ to be the $\bQ(q)$-algebra generated by two elements $i$ and $j$ subject to the relations
$$i^3j-(q^2+q+1)i^2ji+q(q^2+q+1)iji^2-q^3ji^3=0$$
and
$$ij^3-(q^2+q+1)jij^2+q(q^2+q+1)j^2ij-q^3j^3i=0,$$
and set $\CA_q(K)$ to be the $\bQ[q]$-subalgebra generated by $i$ and $j$.

Then $\CA(kK)$ is isomorphic to the specialisation to $q=|k|$ of $\CA_q(K)$.
\end{Thm}

We therefore study the generic composition algebra $\CA_q(K)$.

It is shown in \cite{Szanto} that $\CA_q(K)$ contains elements $\cP(m),\cI(n)$ and $\cR(s)$ for $m,n\geq 0$ and $s\geq1$, with $\cP(0)=j$ and $\cI(0)=i$, satisfying the relations
\begin{enumerate}
\item[1)]
$\begin{aligned}[t]
\cP(m+a)\cP(m) = q^{a+1}&\cP(m)\cP(m+a)\\
&\quad\smash{+(q^{a+1}-q^{a-1})\sum_{r=1}^{\lfloor a/2 \rfloor}\cP(m+r)\cP(m+a-r)}
\end{aligned}$\\
for all $a,m\geq0$;
\item[2)]
$\begin{aligned}[t]
\cI(n)\cI(n+a) = q^{a+1}&\cI(n+a)\cI(n)\\
&\quad\smash{+(q^{a+1}-q^{a-1})\sum_{r=1}^{\lfloor a/2 \rfloor}\cI(n+a-r)\cI(n+r)}
\end{aligned}$\\
for all $a,n\geq0$;
\item[3)]
$\cI(n)\cP(m)=\cR(m+n+1)+q^{m+n}\cP(m)\cI(n)$;
\item[4)]
$\begin{aligned}[t]
\cR(s)\cP(m) = q^s\cP(m)\cR(s) &+ \smash[b]{\sum_{r=1}^{s-1}(q^{s+r}-q^{s+r-2})\cP(m+r)\cR(s-r)}\\
&\qquad\qquad +(q^{2s-1}+q^{2s-2})\cP(m+s);
\end{aligned}$
\item[5)]
$\begin{aligned}[t]
\cI(n)\cR(s) = q^s\cR(s)\cI(n) &+ \smash[b]{\sum_{r=1}^{s-1}(q^{s+r}-q^{s+r-2})\cR(s-r)\cI(n+r)}\\
&\qquad\qquad +(q^{2s-1}+q^{2s-2})\cI(n+s);
\end{aligned}$
\item[6)]
$\cR(s_1)\cR(s_2)=\cR(s_2)\cR(s_1)$.
\end{enumerate}

In particular, these relations are all homogeneous with respect to the grading by dimension vector and can be thought of as defining relations for $\CA_q(K)$. Moreover, we obtain a PBW-type basis for $\CA_q(K)$.
\begin{Thm}[Sz\'ant\'o]\label{CAPBW}
The generic composition algebra $\CA_q(K)$ has a PBW basis
$$\cP(m_1)\cdots\cP(m_a)\cR(s_1)\cdots\cR(s_t)\cI(n_1)\cdots\cI(n_b)$$
such that $m_1\leq\cdots\leq m_a$, $n_1\geq\cdots\geq n_b$ and $s_1\geq\cdots\geq s_t$.
\end{Thm}

We can relate this back to $\CA(kK)$ by specialising to $q=|k|$. Then $\cP(m)\mapsto[P(m)]$, $\cI(n)\mapsto[I(n)]$ and
$$\cR(s)\mapsto \sum_{\substack{(x_1,m_1),\ldots,(x_t,m_t)\\ x_i \mathrm{ distinct}\\ \sum m_id_{x_i}=s}} \prod q^{(m_i-1)d_{x_i}}\frac{(q^{d_{x_i}}-1)}{q-1}[R_{m_1x_1}\oplus\cdots\oplus R_{m_tx_t}].$$

Using this PBW basis and the corresponding relations, we see that when we specialise to $q=0$, we again obtain a PBW basis and have relations
\begin{enumerate}
\item[1)] $\cP(m_2)\cP(m_1)=0$ for $m_1<m_2$;
\item[2)] $\cI(n_2)\cI(n_1)=0$ for $n_1>n_2$;
\item[3)] $\cI(n)\cP(m)=\cR(m+n+1)$ if $m+n>0$; $\cI(0)\cP(0)=\cR(1)+\cP(0)\cI(0)$;
\item[4)] $\cR(s)\cP(m)=0$ if $s>1$; $\cR(1)\cP(m)=\cP(m+1)$;
\item[5)] $\cI(n)\cR(s)=0$ if $s>1$; $\cI(n)\cR(1)=\cI(n+1)$;
\item[6)] $\cR(s_1)\cR(s_2)=\cR(s_2)\cR(s_1)$.
\end{enumerate}

In particular, we deduce in turn that $\cR(1)=ij-ji$, $\cP(m)=(ij-ji)^mj$, $\cI(n)=i(ij-ji)^n$ and $\cR(s+1)=i(ij-ji)^sj$.

We now wish to rewrite these relations in a more transparent form. Let us introduce the notation $\delta_r=i^rj^r$.

We first note that $\cR(s)j=0$ for all $s>1$, using 1) when $m_1=0$. Similarly, $iR(s)=0$ for all $s>1$. Thus the formulae for 6) when $s_1=1$ are equivalent to $\delta_1\cR(s)=\cR(s)\delta_1$ for all $s>1$. Also, it is enough to have 1) in the case $m_2=m_1+1$.

\begin{Lem}
We have that $(ij-ji)^mj=\delta_1^mj-\delta_1^{m-1}j\delta_1-\sum_{r=1}^{m-2}\delta_1^rj\cR(m-r)$.
\end{Lem}

\begin{Pf}
We have that
\begin{align*}
(ij-ji)^{m+1}j &= \delta_1(ij-ji)^mj-ji(ij-ji)^mj\\
&= \delta_1\big(\delta_1^mj-\delta_1^{m-1}j\delta_1-\sum_{r=1}^{m-2}\delta_1^rj\cR(m-r)\big)-j\cR(m+1).
\end{align*}
The proof now follows by induction, the case $m=1$ being clear.
\end{Pf}

We can now rewrite the relations of type 1).
\begin{Lem}
The relations $\cP(m+1)\cP(m)=0$ for all $m\geq 0$ are equivalent to $\delta_1^{m+1}j\delta_1^mj=\delta_1^mj\delta_1^{m+1}j$ for all $m\geq0$, which in turn are equivalent to $\delta_1^mj\delta_1^{m+1}j=\delta_1^{2m+1}j^2$ for all $m\geq 0$.
\end{Lem}

\begin{Pf}
We proceed by induction to prove the first equivalence, the case $m=0$ being clear.

Expanding $\cP(m+1)\cP(m)$ using the previous lemma gives
\begin{align*}
\big(\delta_1^{m+1}&j-\delta_1^mj\delta_1-\sum_{r=0}^{m-1}\delta_1^rj\cR(m+1-r)\big)\big(\delta_1^mj-\delta_1^{m-1}j\delta_1-\sum_{s=0}^{m-2}\delta_1^sj\cR(m-s)\big)\\
&= (\delta_1^{m+1}j\delta_1^mj-\delta_1^mj\delta_1^{m+1}j)-\delta_1(\delta_1^mj\delta_1^{m-1}j-\delta_1^{m-1}j\delta_1^mj)\delta_1\\
&\qquad-\sum_{r=0}^{m-1}\delta_1^rj\cR(m+1-r)\big(\delta_1^mj-\delta_1^{m-1}j\delta_1-\sum_{s=0}^{m-2}\delta_1^sj\cR(m-s)\big)\\
&\qquad\qquad-\sum_{s=0}^{m-2}\delta_1^{m-s}\big(\delta_1^{s+1}j\delta_1^sj-\delta_1^sj\delta_1^{s+1}j\big)\cR(m-s)\\
&=\delta_1^{m+1}j\delta_1^mj-\delta_1^mj\delta_1^{m+1}j,
\end{align*}
using that $\delta_1\cR(s)=\cR(s)\delta_1$ and $\cR(s)j=0$ for all $s>1$, and also that $\delta_1^{r+1}j\delta_1^rj=\delta_1^rj\delta_1^{r+1}j$ for all $0\leq r\leq m$.

Now, we have that
$$\delta_1^{m+1}j\delta_1^mj=\delta_1^2\delta_1^{m-1}j\delta_1^mj,$$
from which the second equivalence follows.
\end{Pf}

In a similar manner, we deduce that the relations of type 4) are redundant. For, after expanding $\cP(m)$, each summand of $\cR(s)\cP(m)$ involves a term of the form
$$\cR(s)\delta_1^rj=\delta_1^r\cR(s)j=0.$$

Analogously, we have that the relations of type 5) are redundant and that we can rewrite the relations of type 2) as $i\delta_1^ni\delta_1^{n+1}=i\delta_1^{n+1}i\delta_1^n$.

Next let us rewrite the relations of type 6).
\begin{Lem}
The relations $\cR(s)\cR(t)=\cR(t)\cR(s)$ for all $s,t\geq 1$ are equivalent to $\delta_s\delta_t=\delta_t\delta_s$ for all $s,t\geq 1$.
\end{Lem}

\begin{Pf}
This time the induction is a bit more involved. Our induction hypotheses will be that $\delta_1^mj=\delta_mj$ for all $m\leq s-1$, that $\cR(m)$ can be expressed using products of $\delta_r$ with $r\leq m$ for all $2\leq m\leq s$ and finally that $\delta_m\delta_n=\delta_n\delta_m$ for all $m,n\leq s$.

Note that in the case $s=1$, there is nothing to show.

Clearly, we always have that $i\delta_rj=\delta_{r+1}$ and that $\delta_r\delta_1j=\delta_{r+1}j$, using the relation $jij^2=ij^3$ repeatedly. N.B. This is the relation of type 1) for $m=0$.

Then
$$\delta_1^sj=\delta_1\delta_1^{s-1}j=\delta_1\delta_{s-1}j=\delta_{s-1}\delta_1j=\delta_sj$$
and thus also
\begin{align*}
\cR(s+1) &= i(ij-ji)^sj=i\delta_1^sj-i\delta_1^{s-1}j\delta_1-\sum_{r=0}^{s-2}i\delta_1^rj\cR(s-r)\\
&=\delta_{s+1}-\delta_s\delta_1-\sum_{r=0}^{s-2}\delta_{r+1}\cR(s-r).
\end{align*}
These give us the first two of our formulae for $s+1$.

It is now clear that the relations $\cR(r)\cR(s+1)=\cR(s+1)\cR(r)$ are equivalent to $\cR(r)\delta_{s+1}=\delta_{s+1}\cR(r)$ for all $1\leq r\leq s$. By considering each $r$ in turn, starting with $r=1$, we see that these relations are equivalent to $\delta_r\delta_{s+1}=\delta_{s+1}\delta_r$ for all $1\leq r\leq s$.
\end{Pf}

\begin{Thm}\label{CArels}
The specialisation $\CA_0(K)$ has defining relations
\begin{enumerate}
\item[1)] $i^mj^{m+1}i^{m+1}j^{m+2}=i^{2m+1}j^{2m+3}$ for all $m\geq 0$;
\item[2)] $i^{n+2}j^{n+1}i^{n+1}j^n=i^{2n+3}j^{2n+1}$ for all $n\geq 0$;
\item[3)] $i^mj^m$ and $i^nj^n$ commute for all $m,n\geq 1$.
\end{enumerate}
Moreover we have a PBW basis in $\CA_0(K)$ involving the elements
$$\cR(1)=ij-ji, \quad \cR(s+1)=i\cR(1)^sj, \quad \cP(m)=\cR(1)^mj, \quad \cI(n)=i\cR(1)^n,$$
as in Theorem \ref{CAPBW}.
\end{Thm}
N.B. These relations again respect the grading by dimension vector.

\begin{Pf}
Let $\mathcal A$ be the $\bQ$-algebra defined by the above relations. We have shown that all the relations hold in $\CA_0(K)$ and thus we have an epimorphism $\mathcal A\to\CA_0(K)$. Moreover, all elements in both $\mathcal A$ and $\CA_0(K)$ can be put into a normal form as in Theorem \ref{CAPBW} with respect to the elements $\cR(s)$, $\cP(m)$ and $\cI(n)$ defined above. Thus for each graded part we have
$$\dim_\bQ\mathcal A_\alpha\leq\dim_{\bQ(q)}\bQ(q)\otimes\CA_q(K)_\alpha=\dim_\bQ\CA_0(K)_\alpha.$$
Thus the map $\mathcal A\to\CA_0(K)$ must be an isomorphism.

N.B. The equality $\dim_{\bQ(q)}\bQ(q)\otimes\CA_q(K)_\alpha=\dim_\bQ\CA_0(K)_\alpha$ can be seen as follows. Each graded part $\CA_q(K)_\alpha$ is a finitely generated free $\bQ[q]$-module. Therefore the dimension when tensored with the quotient field $\bQ(q)$ equals the dimension when specialised to $q=0$, since both equal the rank.
\end{Pf}

The elements $\cR(s)$, $\cP(m)$ and $\cI(n)$ are not well suited to the above description of $\CA_0(K)$. Instead, let us consider the map from the positive root lattice for $K$ to $\CA(K)$ given by sending $\alpha$ to $(\alpha):=i^{\alpha_i}j^{\alpha_j}$.
\begin{Cor}\label{Cor}
The algebra $\CA_q(K)$ has a PBW basis
$$(\dimv P(m_1))\cdots(\dimv P(m_a))\delta_{s_1}\cdots\delta_{s_t}(\dimv I(n_1))\cdots(\dimv I(n_b))$$
such that $m_1\leq\cdots\leq m_a$, $n_1\geq\cdots\geq n_b$ and $s_1\geq\cdots\geq s_t$.
\end{Cor}

\begin{Pf}
Let us first consider the specialisation to $q=0$. The relations defining $\CA_0(K)$ imply that
\begin{align*}
(\dimv P(m_2))(\dimv P(m_1)) &= i^{m_1+m_2}j^{m_1+m_2+2}\\
&= (\dimv P(\lfloor (m_1+m_2)/2\rfloor))(\dimv P(\lceil(m_1+m_2)/2\rceil))
\end{align*}
for all $m_2>m_1$ (and analogously in the preinjective case).

Also, we have that
$$\delta_s(\dimv P(m))=(\dimv P(m+s)) \quad\textrm{and}\quad (\dimv I(n))\delta_s=(\dimv I(n+s))$$
and that
$$(\dimv I(n))(\dimv P(m))=\delta_{n+m+1} \quad\textrm{and}\quad \delta_{s_1}\delta_{s_2}=\delta_{s_2}\delta_{s_1}.$$
Therefore we can put any monomial into the required form. By dimension arguments (working with each graded component separately), we see that this form is unique.

It follows immediately that this basis lifts to the generic composition algebra $\CA_q(K)$.
\end{Pf}

\section{The Composition Monoid}

We now recall the definition of Reineke's composition monoid $\CM(K)$ \cite{Rein2}. In this section, we shall work over an algebraically closed field $k$ of characteristic zero.

Let us write $\Rep(\alpha)$ for the affine space
$$\Rep(\alpha)=\mathbb M(\alpha_j\times\alpha_i)\oplus \mathbb M(\alpha_j\times\alpha_i)$$
and $\GL(\alpha)$ for the affine algebraic group
$$\GL(\alpha_i)\times\GL(\alpha_j).$$
The variety $\Rep(\alpha)$ parametrises the representations of $K$ of dimension vector $\alpha$, the group $\GL(\alpha)$ acts by conjugation on the pairs of matrices and the orbits correspond bijectively to the isomorphism classes of representations of dimension vector $\alpha$.

The extension monoid $\mathcal M(kK)$ is defined to be the free abelian group on irreducible closed $\GL(\alpha)$-stable subvarieties of $\Rep(\alpha)$ for each dimension vector $\alpha$, and with subvarieties $\mathcal A\subset\Rep(\alpha)$ and $\mathcal B\subset\Rep(\beta)$ multiplying according to the rule
$$\mathcal A\ast\mathcal B:=\Big\{X\in\Rep(\alpha+\beta):\begin{gathered}\exists\textrm{ an exact sequence } 0\to B\to X\to A\to 0\\ \textrm{ for some }A\in\mathcal A, B\in\mathcal B\end{gathered}\Big\}.$$

Again, this object seems to be unmanageable and so, by analogy with the Hall algebra, Reineke defines the composition monoid to be the submonoid $\CM(K)$ generated by the varieties $i:=\Rep((1,0))$ and $j:=\Rep((0,1))$.

We shall often abbreviate notion and write $(\alpha)$ for $\Rep(\alpha)=i^{\alpha_i}j^{\alpha_j}$.

The following theorem is proved in \cite{Rein2}.
\begin{Thm}[Reineke]\label{ext}
If $\mathrm{ext}(\mathcal A,\mathcal B)=0$ or if $\mathcal A=\Rep(\alpha)$ and $\mathcal B=\Rep(\beta)$, then
$$\mathrm{codim}\,\mathcal A\ast\mathcal B=\mathrm{codim}\,\mathcal A+\mathrm{codim}\,\mathcal B+\mathrm{ext}(\mathcal B,\mathcal A).$$
\end{Thm}
Here $\mathrm{ext}(\mathcal A,\mathcal B)$ denotes the general (or minimal) value of $\dim\Ext^1(A,B)$ where $A$ and $B$ range over $\mathcal A$ and $\mathcal B$ respectively.

In particular, for the isotropic root $\delta=(1,1)$, we have that $\mathrm{ext}(\delta,\delta)=0$. Thus $\mathrm{ext}(m\delta,n\delta)=0$ for all $m$ and $n$, and thus we have the relation (Corollary 5.5 of \cite{Rein2})
$$\Rep(m\delta)\ast\Rep(n\delta)=\Rep((m+n)\delta),$$
or in other words, $\delta_m\ast\delta_n=\delta_{m+n}$, using the notation $\delta_m=i^mj^m$ as before.

We also have the corollary
\begin{Cor}\label{orbits}
Let $M$ and $N$ be two modules such that $\dim\Ext^1(M,N)=0$. Then
$$\overline{\mathcal O_M}\ast\overline{\mathcal O_N}=\overline{\mathcal O_{M\oplus N}}.$$
\end{Cor}

\begin{Pf}
For an orbit $\mathcal O_M$ in $\Rep(\alpha)$, we have
\begin{align*}
\mathrm{codim}\,\overline{\mathcal O_M} &= \dim\Rep(\alpha)-\dim\GL(\alpha)+\dim\End(M)\\
&= -\langle\alpha,\alpha\rangle+\dim\End(M) = \dim\Ext^1(M,M).
\end{align*}
Here we have written $\langle-,-\rangle$ for the Euler form, given by the matrix $\big(\begin{smallmatrix}1&-2\\0&1\end{smallmatrix}\big)$.

Now, using Theorem \ref{ext}, we have that $\mathrm{codim}\,\overline{\mathcal O_{M\oplus N}}-\mathrm{codim}\,\overline{\mathcal O_M}\ast\overline{\mathcal O_N}$ equals
$$\dim\Ext^1(M\oplus N,M\oplus N)-\dim\Ext^1(M,M)-\dim\Ext^1(N,N)-\dim\Ext^1(N,M),$$
which by assumption equals 0.
\end{Pf}

\begin{Prop}\label{CMrels'}
We have the following relations.
\begin{enumerate}
\item[1)] $(\dimv P(m+a))\ast(\dimv P(m))=(\dimv P(m+\lfloor a/2\rfloor)\ast(\dimv P(m+\lceil a/2\rceil)$ for all $a,m\geq0$;
\item[2)] $(\dimv I(n))\ast(\dimv I(n+a))=(\dimv I(n+\lceil a/2\rceil)\ast(\dimv I(n+\lfloor a/2\rfloor)$ for all $a,n\geq0$;
\item[3)] $(\dimv I(n))\ast(\dimv P(m))=\delta_{m+n}$;
\item[4)] $\delta\ast(\dimv P(m))=(\dimv P(m+1))$;
\item[5)] $(\dimv I(n))\ast\delta=(\dimv I(n+1))$;
\item[6)] $\delta_1^m=\delta_m$.
\end{enumerate}
\end{Prop}

\begin{Pf}
We know that $(\dimv P(m))=\delta_1^mj$, from which the relations of type 4) follow immediately.

Next consider $(\dimv P(m_2))\ast(\dimv P(m_1))$. For all $m_1\leq m_2+1$ we have that $\Ext^1(P(m_1),P(m_2))=0$. Thus Theorem \ref{ext} applies to give
$$(\dimv P(m_2))\ast(\dimv P(m_1))=\Rep(m_1+m_2,m_1+m_2+2).$$
In particular, we have the relations of type 1).

Analogous arguments work for the relations of types 2) and 5).

As mentioned after Theorem \ref{ext}, we have that $\delta_m\ast\delta_n=\delta_{m+n}$ for all $m$ and $n$, and hence we obtain the relations of type 6). The relations 3) now follow immediately.
\end{Pf}

We now prove a normal form for the elements of $\CM(K)$.
\begin{Thm}\label{CMPBW}
Every element of the composition monoid $\CM(K)$ can be expressed uniquely as a product of the form
$$(\dimv P(m_1))\ast\cdots\ast(\dimv P(m_a))\ast\delta_s\ast(\dimv I(n_1))\ast\cdots\ast(\dimv I(n_b))$$
for $m_1\leq\cdots\leq m_a$, $n_1\geq\cdots\geq n_b$ and some $s\geq0$.
\end{Thm}

\begin{Pf}
By definition, $i=(\dimv I(0))$ and $j=(\dimv P(0))$. Thus, using the above relations, we can put any word into the required form. Therefore we just need to prove that distinct words correspond to distinct varieties.

For $m_1\leq\cdots\leq m_a$, consider the product $\overline{\mathcal O_{P(m_1)}}\ast\cdots\ast\overline{\mathcal O_{P(m_a)}}$. By repeated use of Corollary \ref{orbits}, we see that this is precisely $\overline{\mathcal O_{P(m_1)\oplus\cdots\oplus P(m_a)}}$. Analogously for preinjective modules.

That is, we can describe each variety in the composition monoid as
$$\overline{\mathcal O_P}\ast\Rep(s\delta)\ast\overline{\mathcal O_I}$$
for some preprojective module $P$, preinjective module $I$ and some $s\geq 0$.

Now, we know that $\Rep(s\delta)$ contains a dense subset of representations of the form $R=R_{x_1}\oplus\cdots\oplus R_{x_s}$ for distinct $x_1,\ldots,x_s\in\mathbb P^1$. Therefore
$$\mathrm{ext}(\overline{\mathcal O_P},\Rep(s\delta))=0$$
whereas
$$\mathrm{ext}(\Rep(s\delta),\overline{\mathcal O_P})=\dim\Ext(R,P)=s\#\{\textrm{indecomposable summands of $P$}\}.$$
We also have analogous results involving $\Rep(s\delta)$ and $\overline{\mathcal O_I}$. Therefore we can apply Theorem \ref{ext} to get
\begin{align*}
\mathrm{codim}\,\overline{\mathcal O_P}\ast\Rep(s\delta)\ast\overline{\mathcal O_I} &= \dim\Ext(P\oplus R\oplus I,P\oplus R\oplus I)-\dim\Ext(R,R)\\
&= \dim\Ext(P\oplus R\oplus I,P\oplus R\oplus I)-s.
\end{align*}

We shall now show that $\overline{\mathcal O_P}\ast\Rep(s\delta)\ast\overline{\mathcal O_I}$ contains a dense subset of representations isomorphic to $P\oplus R_{x_1}\oplus\cdots\oplus R_{x_s}\oplus I$ for distinct $x_1,\ldots,x_s\in\mathbb P^1$.

Let us fix the modules $P$ and $I$, say corresponding respectively to the pairs of matrices $(P_1,P_2)$ and $(I_1,I_2)$, and set $\alpha:=\dimv P+\dimv I$. For a fixed $s\geq0$ we define a morphism of varieties
$$\phi:\GL(\alpha+s\delta)\times\big(\Rep(\delta)-0\big)^s \to \Rep(\alpha+s\delta)$$
by sending the element $\big(g,(\lambda_1,\mu_1),\cdots,(\lambda_s,\mu_s)\big)$ to $g\cdot(A,B)$, where $A$ and $B$ are the block-diagonal matrices
$$A=\mathrm{diag}(P_1,\lambda_1,\ldots,\lambda_s,I_1) \quad\textrm{and}\quad B=\mathrm{diag}(P_2,\mu_1,\ldots,\mu_s,I_2).$$

We wish to study the fibre over the point $X=g\cdot(A,B)$, and in particular its dimension. So, consider a point $\big(h,(\lambda'_r,\mu'_r)_r\big)$ mapping to $h\cdot(A',B')=X$. For each $r$ let $R_r$ and $R'_r$ be the modules determined by the points $(\lambda_r,\mu_r)$ and $(\lambda'_r,\mu'_r)$ respectively. By the Krull-Remak-Schmidt Theorem we know that there exists a permutation $\pi$ such that $R'_{\pi(r)}\cong R_r$ for each $r$. Thus we can decompose the fibre into a finite number of varieties, one for each permutation, and hence the dimension of the fibre will be the dimension of any of these smaller sets (since they are clearly isomorphic). Also, in the description of the point $X$ we may assume that $g=1$.

That is, we have that $R'_r\cong R_r$ for each $r$ and that $h\cdot(A',B')=(A,B)$.

Now, to say that $R'_r\cong R_r$ is to say that there exists a non-zero $t_r\in k^*$ such that $(t_r\lambda'_r,t_r\mu'_r)=(\lambda_r,\mu_r)$. Thus if we let $T\in\GL(\alpha+s\delta)$ be the point $\big(1,\mathrm{diag}(1,t_1,\ldots,t_s,1)\big)$, then $T\cdot(A',B')=(A,B)$. Therefore, the part of the fibre over $X$ corresponding to the identity permutation is isomorphic to
$$\{(h,t_1,\ldots,t_s):h\in \mathrm{Stab}(A,B)T\}.$$
This is clearly isomorphic to $\mathrm{Stab}(A,B)\times(k^*)^s$, and hence has dimension $\dim\End(P\oplus R\oplus I)+s$, where $R$ is the regular module $R_1\oplus\cdots\oplus R_s$. In particular, the general dimension of a fibre is obtained precisely when the $R_r$ are pairwise non-isomorphic, or equivalently the $x_r:=[\lambda_r:\mu_r]$ are distinct points in $\mathbb P^1$.

We can now calculate the codimension of the closure of the image of $\phi$ using Chevalley's Theorem (see for example \cite{Mumf}).
\begin{align*}
\mathrm{codim}\,&\overline{\mathrm{Im}\phi}\\
&= \dim\Rep(\alpha+s\delta)-\dim\GL(\alpha+s\delta)-2s+\dim\End(P\oplus R\oplus I)+s\\
&= \dim\End(P\oplus R\oplus I)-\langle\alpha+s\delta,\alpha+s\delta\rangle-s\\
&= \dim\Ext(P\oplus R\oplus I)-s.
\end{align*}

Since the image of $\phi$ is contained in the irreducible variety $\overline{\mathcal O_P}\ast\Rep(s\delta)\ast\overline{\mathcal O_I}$ and they have the same codimension, we deduce that they are equal.

Now, by definition, the closure of the image of $\phi$ is precisely the closure of the set
$$\bigcup_{x_1,\ldots,x_s\in\mathbb P^1}\mathcal O_{P\oplus R_{x_1}\oplus\cdots\oplus R_{x_s}\oplus I}.$$
That is, $\overline{\mathcal O_P}\ast\Rep(s\delta)\ast\overline{\mathcal O_I}$ contains a dense set of modules isomorphic to $P\oplus R_{x_1}\oplus\cdots\oplus R_{x_s}\oplus I$ for distinct points $x_1,\ldots,x_s\in\mathbb P^1$.

It follows that distinct expressions of the form $\overline{\mathcal O_P}\ast\Rep(s\delta)\ast\overline{\mathcal O_I}$ give rise to distinct elements in the composition monoid.
\end{Pf}

\begin{Thm}\label{CMrels}
The composition monoid $\CM(K)$ has defining relations
\begin{enumerate}
\item[1)] $i^mj^{m+1}i^{m+1}j^{m+2}=i^{2m+1}j^{2m+3}$ for all $m\geq 0$;
\item[2)] $i^{n+2}j^{n+1}i^{n+1}j^n=i^{2n+3}j^{2n+1}$ for all $n\geq 0$;
\item[3)] $(ij)^m=i^mj^m$ for all $m\geq 1$.
\end{enumerate}
\end{Thm}

\begin{Pf}
It follows from the previous theorem that the relations in Proposition \ref{CMrels'} are defining. Also, it is clear from relation 6) of the proposition that the third, fourth and fifth relations are superfluous and that the first and second are only needed in the cases $a=1$. This leaves precisely the relations written above.
\end{Pf}

\section{Concluding Remarks}

We begin with an immediate corollary of Theorems \ref{CArels} and \ref{CMrels}.
\begin{Cor}\label{epi}
There is a natural surjection from $\CA_0(K)$ to the monoid ring $\bQ\CM(K)$ of the composition monoid. This sends $(\alpha)$ to $\Rep(\alpha)$ for each root $\alpha$. The kernel is generated by the relations
$$(ij)^m=i^mj^m \quad\textrm{for all }m\geq1.$$
\end{Cor}

We note that Theorem \ref{CArels} extends Lemma 4.2 in \cite{Rein2}. That is, only the quantum Serre relations are considered in \cite{Rein2}, whereas we have exhibited defining relations for the specialisation to $q=0$ of the composition algebra. Similarly, the above corollary extends Proposition 4.3 of \cite{Rein2}.

We stress that, whereas in the composition algebra the subalgebra generated by the elements $\delta_m$ form a polynomial algebra with infinitely many variables (one for each $m$), the corresponding subalgebra in the monoid ring $\bQ\CM(K)$ is isomorphic to a polynomial algebra on the single variable $\delta$.

At the end of \cite{Rein2}, Reineke poses four questions concerning the composition monoid. Using our results, we can answer all of these for the case of the Kronecker quiver.

The first two questions concern defining relations for the composition monoid. We have given in Theorem \ref{CMrels}, however, a set of minimal defining relations, which answers the second question. These relations can be written as
$$\Rep(\dimv P(m))\ast\Rep(\dimv P(m+1))=\Rep(\dimv P(m)\oplus P(m+1)) \textrm{ for all }m\geq 0,$$
and similarly for the prinjective modules, as well as
$$\Rep(\delta)\ast\Rep(m\delta)=\Rep((m+1)\delta) \quad\textrm{for all }m\geq 1.$$
Therefore we have an affirmative answer to Reineke's first question --- namely that only relations of this type are necessary.

The third question discusses a partial normal form for elements of the composition monoid. It is proved in Theorem 5.8 of \cite{Rein2} that every word can be expressed in the form
$$\Rep(\alpha_1)\ast\cdots\ast\Rep(\alpha_a),$$
where each $\alpha_r$ is a Schur root and such that
$$\mathrm{ext}(\alpha_r,\alpha_{r+1})\neq0 \quad\textrm{implies}\quad \mathrm{ext}(\alpha_{r+1},\alpha_r)\neq0.$$

We recall that $\alpha$ is a Schur root if there exists a representation $X$ of dimension vector $\alpha$ with $\End(X)=k$.

Now, the normal form that we have described in Theorem \ref{CMPBW} is precisely of this type. In fact, for the Kronecker quiver, we always have that $\mathrm{ext}(\alpha_r,\alpha_{r+1})=0$. Moreover, as was shown in the proof of Theorem \ref{CMPBW}, we can describe the elements of $\CM(K)$ just in terms of Schur roots.

The Schur roots for the Kronecker quiver are precisely the dimension vectors of the indecomposable preprojectives and preinjectives together with $\delta$. We order these by saying that
$$\dimv P(m)<\dimv P(m+1)<\delta<\dimv I(n+1)<\dimv I(n) \quad\textrm{for all }m,n\geq0.$$
Note that $\mathrm{ext}\,(\alpha,\beta)=0$ for Schur roots $\alpha<\beta$.

The elements of $\CM(K)$ now correspond precisely to finite sets of Schur roots. In particular, $\alpha_1\leq\cdots\leq\alpha_a$ corresponds to the variety
$$\mathcal A:=\Rep(\alpha_1)\ast\cdots\ast\Rep(\alpha_a)$$
and the general element of $\mathcal A$ is of the form $R=R_1\oplus\cdots\oplus R_a$ with $\dimv R_r=\alpha_r$. The codimension of $\mathcal A$ can also be calculated using Theorem \ref{ext} to be
$$\mathrm{codim}\,\mathcal A=\sum_{s>r}\mathrm{ext}(\alpha_s,\alpha_r).$$

In this way we recover the Kac's generic decomposition of $\alpha$ \cite{Kac1} (see also \cite{Scho}) as the decomposition of $\Rep(\alpha)$. That is, the decomposition
$$\alpha=\alpha_1+\cdots+\alpha_a \quad\textrm{with each $\alpha_r$ a Schur root and }\mathrm{ext}(\alpha_r,\alpha_s)=0\textrm{ for }r\neq s.$$

Finally we mention the fourth question raised by Reineke, which asks whether the partial normal form can be lifted to the composition algebra.

Let us order the positive roots of $K$ as follows.
$$\dimv P(m)<\dimv P(m+1)< (s+1)\delta<s\delta< \dimv I(n+1)<\dimv I(n)$$
for all $m,n\geq 0$ and $s\geq1$. Corollary \ref{Cor} then states that the elements
$$(\alpha_1)\cdots(\alpha_a) \quad\textrm{for roots }\alpha_1\leq\cdots\leq\alpha_a$$
form a PBW basis for $\CA_0(K)$, where $(\alpha)=i^{\alpha_i}j^{\alpha_j}$. It is now immediate that the image in the composition monoid of each of these basis elements is already in normal form.

Moreover, we see that we have distinct basis elements $i^mj^m$ for each $m\geq 1$. These arise since there exist indecomposable representations of dimension vector $m\delta$ over the finite field $\mathbb F_{\!q}$ with endomorphism ring $\mathbb F_{\!q^m}$. This phenomenon only occurs for non-algebraically closed fields, and so does not appear in the composition monoid. This explains why the epimorphism described in Corollary \ref{epi} is not an isomorphism.

\end{document}